\documentclass{article}
\usepackage[centertags]{amsmath}

\begin{document}

\title{Counterexamples on Jumarie's three basic fractional calculus formulae for non-differentiable continuous functions}

\author{Cheng-shi Liu\\Department of Mathematics\\Northeast Petroleum University\\Daqing 163318, China
\\Email: chengshiliu-68@126.com}

 \maketitle

\begin{abstract}
  Juamrie proposed a modified Riemann-Liouville derivative definition and gave three so-called basic fractional calculus
formulae
  $(u(t)v(t))^{(\alpha)}=u^{(\alpha)}(t)v(t)+u(t)v^{(\alpha)}(t)$, $(f(u(t)))^{(\alpha)}=f'_uu^{(\alpha)}(t)$ and $(f(u(t)))^{(\alpha)}=(f(u))^{(\alpha)}(u'(t))^{\alpha}$ where $u$ and $v$ are required to be non-differentiable and continuous for the first formula, $f$ is assumed to be differentiable for the second formula, while in the third formula $f$ is non-differentiable and $u$ is differentiable, at the point $t$.
  I once gave three counterexamples to show that Jumarie's three formulae are not true for differentiable functions(Cheng-shi Liu. Counterexamples on Jumarie¡¯s two basic fractional calculus formulae. Communications in Nonlinear Science and Numerical Simulation, 2015, 22(1): 92-94.). However, these examples cannot directly become the suitable counterexamples for the case of non-differentiable continuous functions. In the present paper, I first provide five counterexamples to show directly  the Jumarie's formulae are also not true for non-differentiable continuous functions. Then I prove that essentially non-differentiable cases can be transformed to the differentiable cases.  Therefore, those counterexamples in the above paper are indirectly right. In summary, the Jumarie's formulae are not true. This paper can be considered as the corrigendum and supplement to the above paper.


\end{abstract}

\section{Introduction}

 Jumarie proposed a modified Riemann-Liouville fractional
 derivative[1-5]:
 \begin{equation}
  f^{(\alpha)}(t)=\frac{1}{\Gamma(1-\alpha)}\frac{\mathrm{d}}{\mathrm{d}t}
  \int_0^t(t-x)^{-\alpha}(f(x)-f(0))\mathrm{d}x.
 \end{equation}
and gave some basic fractional calculus formulae(see, for
example, formulae (3.11)-(3.13) in Ref.[4] or formulae (4.3), (4.4) and (4.5)in Ref.[5]):
 \begin{equation}
 (u(t)v(t))^{(\alpha)}=u^{(\alpha)}(t)v(t)+u(t)v^{(\alpha)}(t),
 \end{equation}
\begin{equation}
 (f(u(t)))^{(\alpha)}=f'_uu^{(\alpha)}(t),
 \end{equation}
where Jumarie requires the functions $u$ and $v$ are non-differentiable and continuous, while $f$ is differentiable, at the point $t$. Jumarie's third formula is given by
\begin{equation}
(f(u(t)))^{(\alpha)}=(f(u))^{(\alpha)}(u'(t))^{\alpha},
\end{equation}
where $f$ is non-differentiable and $u$ is differentiable at the point $t$.

The formula (3) has been applied to solve the exact solutions to some
 nonlinear fractional order differential equations(see, for example, Refs.[6-9]).

 In [10], I once gave three counterexamples to show that  Jumarie's these  so-called basic
  formulae are not correct. But, I neglected the conditions of the Jumarie's formulae. Indeed, for example, the formula (2) requires that the functions $u$ and $v$ are non-differentiable and continuous, and the formula (3) requires that $f$ is differentiable while $u$ is non-differentiable and continuous, at the point $t$. Therefore, the examples in [10] cannot be considered as suitable direct counterexamples to Jumarie's formulae under the conditions of non-differentiable functions. However, I will show that these counterexamples do hold indirectly. In addition, since the functions in [6-9] need to be differentiable, my counterexamples are still right for these applications. In [11], Jumarie emphasizes that it is just at some point that his formulae do hold. At such point, the function is non-differentiable and continuous, and the fractional derivative exists. In the present paper, I further provide five counterexamples which satisfy all conditions in Jumarie's formulae to show directly that Jumarie's formulae are incorrect for the cases of non-differentiable continuous functions. Finally, I prove that essentially non-differentiable cases can be transformed to the differentiable cases. Therefore, in other words, those counterexamples in the paper [10] are indirectly valid.

Recently, some problems about the rules of fractional derivatives have been discussed by some authors(see, for example, [10-13]). For instance, Tarasov[12] gave an important result for Leibniz rule. For local fractional derivative of nowhere differentiable continuous functions on open intervals, some detailed discussions can be found in [13].

 \section{Counterexamples to formula (2)}
As in [10], we need the $\frac{1}{2}-$order derivatives of the following four
functions $f(t)=t, f(t)=\sqrt t$, $f(t)=t^2$ and $f(t)=t^\frac{3}{2}$ with $f(0)=0$:
\begin{equation}
(t)^{(1/2)}=2\sqrt{\frac{t}{\pi}},
\end{equation}
\begin{equation}
(\sqrt{t})^{(1/2)}=\frac{\sqrt{\pi }}{2},
\end{equation}
\begin{equation}
(t^2)^{(1/2)}=\frac{8t^{3/2}}{3\sqrt{\pi}},
\end{equation}
\begin{equation}
(t^\frac{3}{2})^{(1/2)}=\frac{3\sqrt{\pi}t}{4}.
\end{equation}

 \textbf{ Counterexample 1 }(The counterexample of formula (2)). Take  $\alpha=\frac{1}{2}$ and
 \begin{equation}
    u(t)=
   \begin{cases}
   \sqrt{t}, &\mbox{$0\leq t\leq1$,}\\
  \sqrt{t}+t-1, &\mbox{$t>1$}.
   \end{cases}
  \end{equation}
It is easy to see that $u(t)$ is continuous, and is non-differentiable at $t=1$. Further, we have
\begin{equation}
    H(t)= \int_0^t(t-x)^{-\alpha}(f(x)-f(0))\mathrm{d}x=
   \begin{cases}
  \int_0^t\frac{\sqrt x}{\sqrt{t-x}}\mathrm{d}x, &\mbox{$0\leq t\leq1$,}\\
  \int_0^1\frac{\sqrt x}{\sqrt{t-x}}\mathrm{d}x+ \int_1^t\frac{\sqrt x+x-1}{\sqrt{t-x}}\mathrm{d}x, &\mbox{$t>1$}.
   \end{cases}
  \end{equation}
And then, we have
\begin{equation}
    H(t)=
   \begin{cases}
  \int_0^t\frac{\sqrt x}{\sqrt{t-x}}\mathrm{d}x, &\mbox{$0\leq t\leq1$,}\\
  \int_0^t\frac{\sqrt x}{\sqrt{t-x}}\mathrm{d}x+ \int_1^t\frac{x-1}{\sqrt{t-x}}\mathrm{d}x, &\mbox{$t>1$}.
   \end{cases}
  \end{equation}
Let $t-x=s^2$. Then we have
\begin{equation}
K(t)=\int_1^t\frac{x-1}{\sqrt{t-x}}\mathrm{d}x=2\int_0^{\sqrt{t-1}}(t-1-s^2)\mathrm{d}x=\frac{4}{3}(t-1)^{\frac{3}{2}}.
\end{equation}
Therefore, if $0\leq t<1$,
\begin{equation}
 u^{(1/2)}(t)=(\sqrt t)^{(\frac{1}{2})}=\frac{\sqrt{\pi}}{2};
 \end{equation}
If $t>1$,
 \begin{equation}
 u^{(1/2)}(t)=(\sqrt t)^{(\frac{1}{2})}+\frac{1}{\sqrt\pi}K'(t)=\frac{\sqrt{\pi }}{2}+\frac{2}{\sqrt\pi}(t-1)^{\frac{1}{2}},
 \end{equation}
where we use $\Gamma(\frac{1}{2})=\sqrt\pi$. And hence, at $t=1$, it follows that $ u^{(1/2)}(1)$ exists and
\begin{equation}
 u^{(1/2)}(1)=\frac{\sqrt{\pi}}{2}.
 \end{equation}
Further, by taking $v(t)=u(t)$, we have
 \begin{equation}
 u^{(1/2)}(1)v(1)+u(1)v^{(1/2)}(1)=\sqrt{\pi}.
 \end{equation}
On the other hand, we have
\begin{equation}
    u(t)v(t)=
   \begin{cases}
   t, &\mbox{$0\leq t\leq1$,}\\
  (\sqrt{t}+t-1)^2, &\mbox{$t>1$}.
   \end{cases}
  \end{equation}
Hence, if $t<1$, we have
\begin{equation}
 (uv)^{(1/2)}(t)=(t)^{(\frac{1}{2})}=2\sqrt{\frac{t}{\pi}};
 \end{equation}
If $t>1$, we have
 \begin{equation}
 (uv)^{(1/2)}(t)= \frac{1}{\sqrt\pi}\frac{\mathrm{d}}{\mathrm{d}t}\{\int_0^1\frac {x}{\sqrt{t-x}}\mathrm{d}x+ \int_1^t\frac{(\sqrt x+x-1)^2}{\sqrt{t-x}}\mathrm{d}x\}.
 \end{equation}
Further, we have
 \begin{equation}
 (uv)^{(1/2)}(t)= (t)^{(\frac{1}{2})}+\frac{1}{\sqrt\pi}\int_1^t\frac{3\sqrt x+2(x-1)-x^{-\frac{1}{2}}}{\sqrt{t-x}}\mathrm{d}x.
 \end{equation}
By computing the last integral, we get
 \begin{equation}
 (uv)^{(1/2)}(t)=2\sqrt{\frac{t}{\pi}}+ \frac{1}{\sqrt\pi}\{\frac{8}{3}(t-1)^{\frac{3}{2}}+3\sqrt{t-1}+3t(\frac{\pi}{2}-\arcsin{\frac{1}{\sqrt t}})+2\arcsin{\frac{1}{\sqrt t}}-\pi\}.
 \end{equation}
Therefore, at $t=1$, $ (uv)^{(1/2)}(t)$ does exist and $ (uv)^{(1/2)}(1)=\frac{2}{\sqrt\pi}\neq\sqrt\pi$. From (16), it turns out that we  have
 \begin{equation}
 (uv)^{(1/2)}(1)\neq u^{(1/2)}(1)v(1)+u(1)v^{(1/2)}(1).
 \end{equation}
This example shows that Jumarie's formula (2) is not true for the non-differentiable continuous functions.

Next, we give a more simple example.

\textbf{ Counterexample 2}. Take $\alpha=\frac{1}{2}$ and
 \begin{equation}
    u(t)=
   \begin{cases}
  1-t, &\mbox{$ t\leq1$,}\\
  t-1, &\mbox{$t>1$}.
   \end{cases}
  \end{equation}
It is easy to see that $u(t)$ is continuous, and is non-differentiable at $t=1$. Further, we have
\begin{equation}
    H(t)= \int_0^t(t-x)^{-\alpha}(u(x)-u(0))\mathrm{d}x=
   \begin{cases}
  \int_0^t\frac{-x}{\sqrt{t-x}}\mathrm{d}x, &\mbox{$t\leq1$,}\\
  \int_0^1\frac{-x}{\sqrt{t-x}}\mathrm{d}x+ \int_1^t\frac{x-2}{\sqrt{t-x}}\mathrm{d}x, &\mbox{$t>1$}.
   \end{cases}
  \end{equation}
 And then, we have
\begin{equation}
    H(t)=
   \begin{cases}
  \int_0^t\frac{-x}{\sqrt{t-x}}\mathrm{d}x, &\mbox{$t\leq1$,}\\
  \int_0^t\frac{-x}{\sqrt{t-x}}\mathrm{d}x+ 2\int_1^t\frac{x-1}{\sqrt{t-x}}\mathrm{d}x, &\mbox{$t>1$}.
   \end{cases}
  \end{equation}
 Therefore, if $ t<1$,
\begin{equation}
 u^{(1/2)}(t)=-(t)^{(\frac{1}{2})}=-2\sqrt{\frac{t}{\pi}};
 \end{equation}
If $t>1$,
 \begin{equation}
 u^{(1/2)}(t)=\frac{1}{\sqrt\pi}H'(t)=-2\sqrt{\frac{t}{\pi }}+4\sqrt{\frac{t-1}{\pi }}.
 \end{equation}
It follows that
 \begin{equation}
 u^{(1/2)}(1)=-\frac{2}{\sqrt\pi}.
 \end{equation}
Hence, from $u(1)=0$ we have
 \begin{equation}
 2u(1)u^{(1/2)}(1)=0.
 \end{equation}
 On the other hand, we have $u^2(t)=(t-1)^2=t^2-2t+1$, and then
 \begin{equation}
 (u^2)^{(1/2)}(t)=(t^2)^{(1/2)}-2(t)^{(1/2)}=\frac{8}{3\sqrt\pi}t^{\frac{3}{2}}-\frac{4}{\sqrt\pi}t^{\frac{1}{2}}.
 \end{equation}
 Therefore, we get
 \begin{equation}
 (u^2)^{(1/2)}(1)=-\frac{4}{3\sqrt\pi}\neq0.
 \end{equation}
 So we give
  \begin{equation}
 (u^2)^{(1/2)}(1)\neq 2u(1)u^{(1/2)}(1).
 \end{equation}
Therefore, if we take $v(t)=u(t)$, we have equivalently from (32)
 \begin{equation}
 (uv)^{(1/2)}(1)\neq u^{(1/2)}(1)v(1)+u(1)v^{(1/2)}(1).
 \end{equation}
This shows again that Jumarie's formula (2) is not true for the non-differentiable continuous functions.

\section{Counterexamples to formula (3)}

\textbf{ Counterexample 3 }(The counterexample of formula (3)). Take
$f(u)=u^2$ and $u(t)$ is also given by (9), and $\alpha=\frac{1}{2}$. Then $f(u(t))$ is given by
\begin{equation}
    f(u(t))=
   \begin{cases}
   t, &\mbox{$0\leq t\leq1$,}\\
  (\sqrt{t}+t-1)^2, &\mbox{$t>1$}.
   \end{cases}
  \end{equation}

 Firstly, therefore, we get the $\frac{1}{2}$ order derivative of $f(u(t))$ at $t=1$
 \begin{equation}
 (f(u(t)))^{(1/2)}|_{t=1}=(u^2(t))^{(1/2)}|_{t=1}=\frac{2}{\sqrt{\pi}}.
 \end{equation}
 Secondly, from (15) and $u(1)=1$, at $t=1$, we have
  \begin{equation}
 f'_u(1)u^{(1/2)}(1)=\sqrt{\pi}.
 \end{equation}
So we find
\begin{equation}
 (f(u(t)))^{(1/2)}|_{t=1}\neq f'_uu^{(1/2)}(t)|_{t=1}.
 \end{equation}
This shows that the Jumarie's  formula (3) is not correct for non-differentiable continuous functions.

\textbf{Counterexample 4}. Take $u(t)$ as the same as (23) in counterexample 2, and $f(u)=u^2$. According to the counterexample 2,  we have
\begin{equation}
 (f(u(t)))^{(1/2)}|_{t=1}=(u^2)^{(1/2)}(1)=-\frac{4}{3\sqrt\pi},
 \end{equation}
 and
 \begin{equation}
 f'(u)u^{(1/2)}(t)|_{t=1}=2u(1)u^{(1/2)}(1)=0.
 \end{equation}
 It turns out that
 \begin{equation}
 (f(u(t)))^{(1/2)}|_{t=1}\neq f'_uu^{(1/2)}(t)|_{t=1}.
 \end{equation}
This shows again that the Jumarie's formula (3) is not correct for non-differentiable continuous functions.

\section{Counterexamples to formula (4)}

\textbf{Counterexample 5}. We take $\alpha=\frac{1}{2}$, $u(t)=t^2$. $f(u)$ is taken as the following form
 \begin{equation}
    f(u)=
   \begin{cases}
   \sqrt{u}, &\mbox{$0\leq u\leq1$,}\\
  \sqrt{u}+u-1, &\mbox{$u>1$}.
   \end{cases}
  \end{equation}
It is easy to see that $f(u)$ is non-differentiable at $u=1$ and $f^{(\frac{1}{2})}(1)=\frac{\sqrt\pi}{2}$ from the formula (15). Therefore, we have
\begin{equation}
f^{(\frac{1}{2})}(u)(u'(t)))^{\frac{1}{2}}|_{t=1}=\sqrt{\frac{\pi}{2}}.
\end{equation}
On the other hand, since $u(t)=t^2$, we have
\begin{equation}
    f(u(t))=
   \begin{cases}
   t, &\mbox{$0\leq t\leq1$,}\\
  t+t^2-1, &\mbox{$t>1$}.
   \end{cases}
  \end{equation}
Therefore, when $0\leq t<1$, we have
\begin{equation}
(f(u(t)))^{(1/2)}=(t)^{(1/2)}=2\sqrt{\frac{t}{\pi}}.
\end{equation}
When $t>1$, we have
\begin{equation}
(f(u(t)))^{(1/2)}=(t)^{(1/2)}=2\sqrt{\frac{t}{\pi}}+\frac{2}{\sqrt\pi}(t\sqrt{t-1}-\frac{1}{3}(t-1)^{\frac{3}{2}}).
\end{equation}
In fact, when $t>1$, we have
\begin{equation}
    H(t)= \int_0^t(t-x)^{-\alpha}(f(u(x))-f(u(0)))\mathrm{d}x=
   \int_0^1\frac{ x}{\sqrt{t-x}}\mathrm{d}x+ \int_1^t\frac{x+x^2-1}{\sqrt{t-x}}\mathrm{d}x.
  \end{equation}
And then, we have
\begin{equation}
    H(t)=
  \int_0^t\frac{x}{\sqrt{t-x}}\mathrm{d}x+ \int_1^t\frac{x^2-1}{\sqrt{t-x}}\mathrm{d}x.
  \end{equation}
hence, we have
\begin{equation}
(f(u(t)))^{(1/2)}=\frac{1}{\sqrt\pi}H'(t)=(t)^{(1/2)}+\frac{1}{\sqrt\pi}\int_1^t\frac{2x}{\sqrt{t-x}}\mathrm{d}x.
\end{equation}
By taking $x=t-s^2$, we have
\begin{equation}
\int_1^t\frac{x}{\sqrt{t-x}}\mathrm{d}x=\int_0^{\sqrt{t-1}}(t-s^2)\mathrm{d}s=t\sqrt{t-1}-\frac{1}{3}(t-1)^{\frac{3}{2}}.
\end{equation}
It follows that
\begin{equation}
(f(u(t)))^{(1/2)}=2\sqrt{\frac{t}{\pi}}+2t\sqrt{t-1}-\frac{2}{3}(t-1)^{\frac{3}{2}}.
\end{equation}
Therefore, we get
\begin{equation}
(f(u(t))^{(1/2)}|_{t=1}=\frac{2}{\sqrt\pi}.
\end{equation}
By (42) and (51), it turns out that
\begin{equation}
(f(u(t)))^{(1/2)}|_{t=1}\neq (f(u))^{(1/2)}(u'(t))^{\frac{1}{2}}|_{t=1}.
\end{equation}
This shows  that the Jumarie's formula (4) is not true for non-differentiable continuous functions.

\section{Explanation}

\textbf{Theorem}. Let $u_1(t)$ and $u_2(t)$ be two continuous differentiable functions on $[0,+\infty)$, that is, belong to $C^1[0,+\infty)$. Take the function $u(t)$ as
\begin{equation}
    u(t)=
   \begin{cases}
  u_1(t), &\mbox{$ t\leq t_0$,}\\
  u_2(t), &\mbox{$t>t_0$},
   \end{cases}
  \end{equation}
where $t_0>0$ is a fixed point,  $u_1(t_0)=u_2(t_0)$ and $u_1'(t_0)\neq u_2'(t_0)$, that is, $u(t)$ is continuous and non-differentiable at the point $t_0$. If $u_1^{(\alpha)}(t)$ is continuous, and $u_1^{(\alpha)}(t_0)$ exists and is finite, then $u^{(\alpha)}(t_0)$ also exists and $u^{(\alpha)}(t_0)=u_1^{(\alpha)}(t_0)$, where $0<\alpha<1$.

\textbf{Proof}. By the definition of $u(t)$ and the conditions of the theorem, we have
\begin{equation}
    u^{(\alpha)}(t)=
   \begin{cases}
  u_1^{(\alpha)}(t), &\mbox{$ t<t_0$,}\\
  \frac{1}{\Gamma(1-\alpha)}\frac{\mathrm{d}}{\mathrm{d}t}\{\int_0^{t_0}\frac{u_1(x)-u_1(0)}{(t-x)^{\alpha}}\mathrm{d}x+
  \int_{t_0}^t\frac{u_2(x)-u_1(0)}{(t-x)^{\alpha}}\mathrm{d}x\}, &\mbox{$t>t_0$}.
   \end{cases}
  \end{equation}
Further, we have
\begin{equation}
    u^{(\alpha)}(t)=
   \begin{cases}
  u_1^{(\alpha)}(t), &\mbox{$ t<t_0$,}\\
  \frac{1}{\Gamma(1-\alpha)}\frac{\mathrm{d}}{\mathrm{d}t}\{\int_0^{t}\frac{u_1(x)-u_1(0)}{(t-x)^{\alpha}}\mathrm{d}x+
  \int_{t_0}^t\frac{u_2(x)-u_1(x)}{(t-x)^{\alpha}}\mathrm{d}x\}, &\mbox{$t>t_0$}.
   \end{cases}
  \end{equation}
 Therefore, we get
 \begin{equation}
    u^{(\alpha)}(t)=
   \begin{cases}
  u_1^{(\alpha)}(t), &\mbox{$ t<t_0$,}\\
   u_1^{(\alpha)}(t)+\frac{1}{\Gamma(1-\alpha)}\frac{\mathrm{d}}{\mathrm{d}t}
  \int_{t_0}^t\frac{u_2(x)-u_1(x)}{(t-x)^{\alpha}}\mathrm{d}x, &\mbox{$t>t_0$}.
   \end{cases}
  \end{equation}
Since $u_2(x)-u_1(x)$ is differentiable, we have
\begin{equation}
\frac{\mathrm{d}}{\mathrm{d}t}
  \int_{t_0}^t\frac{u_2(x)-u_1(x)}{(t-x)^{\alpha}}\mathrm{d}x= \int_{t_0}^t\frac{u_2'(x)-u_1'(x)}{(t-x)^{\alpha}}\mathrm{d}x.
\end{equation}
Because $u_2'(x)-u_1'(x)$ is continuous from the conditions of the theorem, and $\frac{1}{(t-x)^{\alpha}}$ is integrable in $[t_0,t]$, by the mean value theorem of integral(see, for example,[14]), we have
\begin{equation}
  \int_{t_0}^t\frac{u_2'(x)-u_1'(x)}{(t-x)^{\alpha}}\mathrm{d}x=
  (u_2'(\xi)-u_1'(\xi))\frac{(t-t_0)^{1-\alpha}}{1-\alpha},
\end{equation}
where $\xi\in(t_0,t)$.  Hence, we obtain
\begin{equation}
\lim_{t\rightarrow t_0}
 \frac{\mathrm{d}}{\mathrm{d}t}
  \int_{t_0}^t\frac{u_2(x)-u_1(x)}{(t-x)^{\alpha}}\mathrm{d}x=
  \lim_{t\rightarrow t_0}(u_2'(\xi)-u_1'(\xi))\frac{(t-t_0)^{1-\alpha}}{1-\alpha}=0.
\end{equation}
From the continuity of $u_1^{(\alpha)}(t)$ at the point $t_0$, it follows that
\begin{equation}
u^{(\alpha)}(t_0)=u_1^{(\alpha)}(t_0).
\end{equation}
The proof is completed.

According to the above theorem, we know that $u^{(\alpha)}(t_0)$ doesn't depend on the function $u_2$, that is, $u^{(\alpha)}(t_0)$ doesn't depend on the values of $u(t)$ on $(t_0,\infty)$. On the other hand, from the definition (1), we can also see the fact. Therefore, to compute $u^{(\alpha)}(t_0)$, we can smoothly continue the function $u(t)$ from $u_1(t)$ on $[0,t_0]$ to $(t_0,\infty)$, such that $u(t)$ is differentiable at the point $t_0$. In other words, essentially, those counterexamples in [10] are also right. For instance, the function $u(t)=\sqrt t$ in the counterexample 1 in [10] is just the smooth continuation of the function $\sqrt t $ defined on $[0,1]$ in the counterexample 1 in the present paper.

\section{conclusion}

By the above counterexamples and theory, we have showed that the Jumarie's three basic formulae for fractional derivative are not correct for non-differentiable continuous functions.

\textbf{Acknowledgements}. Thanks to Prof.Jumarie and Prof.Kamocki for their pointing out to me the negligence in my paper[12] on the conditions of the Jumarie's formulae.

\end{document}